\documentstyle[11pt]{article} \openup 5pt 
\oddsidemargin
0pt \evensidemargin 0pt \topmargin -30pt \textwidth 6.6truein
\textheight
9truein \parskip 12pt{}\baselineskip 10pt{}\lineskip 10pt{}\def \VSPACE
{{\!\!\!\!}_{_{_{_{_{_{_{_{_{_{_{_{\,}}}}}}}}}}}}}
\def\OPLUS#1#2{\VSPACE
\mbox{\raisebox{-5pt}{${}^{\dis{#1}\oplus} _{#2}$}}} \def\INF{\infty}
 \def \Z{\hbox{$Z\hskip -5.2pt Z$}} \def \LL
{\hbox{$I\hskip -5.2pt L$}} \def\sZ{\hbox{$\sc Z\hskip -4.2pt Z$}}   \def \C{\hbox{$C\hskip -5pt \vrule height 6pt
depth 0pt \hskip 6pt$}} \def \sC{{\hbox{$\sc C\hskip -5pt \vrule height
5pt depth 0pt \hskip 6pt$}}} \def\qed{\hfill \hfill \ifhmode\unskip
\nobreak\fi\ifmmode\ifinner \else\hskip5pt\fi\fi \hbox{\hskip5pt\vrule
width4pt height6pt depth1.5pt\hskip 1 pt}}  
\def\d{\delta} \def\D{\Delta}  
\def\l{\lambda
} \def\L{\Lambda}      \def\Vir{\hbox{\bf\sl Vir}}
\def\HVir{\hbox{{\bf
\sl Vir}[$M$]}} \def\HVira{\hbox{{\bf\sl Vir}[$M_1$]}}
\def\HVB{\hbox{{\bf
\sl$\overline{\mbox{Vir}}$}[$M$]}} \def\sc{\scriptstyle} \def\ssc
{\scriptscriptstyle} \def\dis{\displaystyle} \def\cl{\centerline}
  \def\nl{\newline} \def\ol{\overline}
\def
\ul{\underline} \def\wt{\widetilde} 
\def\rar{\rightarrow
}   
\def\Rla
{\Leftrightarrow}  \def\bs{\backslash} \def\hs
{\hspace*} \def\rb{\raisebox} \def\VS{\mbox{\rb{-10pt}{\,}}}   \def\ni{\noindent} \def\hi{\hangindent}
\def\ha
{\hangafter} \def\iso{\cong} 
\begin{document} \cl{\bf CLASSIFICATION OF
SIMPLE HARISH-CHANDRA MODULES }\cl{\bf OVER THE HIGH RANK VIRASORO
ALGEBRAS }\par\ \par\ \par\cl{ Yucai Su }\par\ \par\cl{ Department of
Applied Mathematics, Shanghai Jiaotong University, China }\cl{ Email:
kfimmi@public1.sta.net.cn }\par\ \par\ \par\ni{\bf ABSTRACT.} A notion
of
generalized highest weight modules over the high rank Virasoro algebras
is introduced in this paper, and a theorem, which was originally given
as
a conjecture by Kac over the Virasoro algebra, is generalized. Mainly,
we
prove that a simple Harish-Chandra module over a high rank Virasoro
algebra is either a generalized highest weight module, or a module of
the
intermediate series. \par\ni{\bf KEYWORDS:} High rank Virasoro algebra,
Harish-Chandra module, uniformly bounded module, generalized highest
weight module, module of the intermediate series.\par\ \par\ \par\cl{\bf
\S1. Introduction}\par The Virasoro algebra \Vir, closely related to
Kac-Moody algebras [{\ul 5},{\ul 6}], is of interest to both
mathematicians and physicists (see, e.g., [{\ul 2},{\ul 3},{\ul 4},{\ul
8},{\ul{12}}]). This is partly due to its relevance to string theories
[\ul{12}] and 2-dimensional conformal field theories [{\ul 3}]. As the
universal central extension of an infinite dimensional complex Lie
algebra of the linear differential operators $\{t^{i+1}{d\over{dt}}\ |\
i
\in\Z\}$, \Vir\, can be realized by taking a basis $\{L_i,c\ |\
i\in\Z\}$
with the following commutation relations:
$$\matrix{[L_i,L_j]=(j-i)L_{i+j}
+{i^3-i\over12}\d_{i,-j}c,\VS\hfill\cr [L_i,c]=[c,c]=0,\hfill\cr}\ \ \ \
\ i,j\in\Z.\eqno(1.1)$$\par For an $n$-dimensional \Z-submodule $M=M_n$
of \C, the rank $n$ Virasoro algebra \HVir\, (also called a high rank
Virasoro algebra when $n\ge 2$), a notion generalized from that of the
Virasoro algebra and introduced in [\ul{11}], is defined as the complex
Lie algebra with basis $\{L_\mu,c\ |\ \mu\in M\}$ and the commutation
relations
$$\matrix{[L_\mu,L_\nu]=(\nu-\mu)L_{\mu+\nu}-{\mu^3-\mu\over12}
\d_{\mu,-\nu}c,\VS\hfill\cr [L_\mu,c]=[c,c]=0,\hfill\cr}\ \ \ \ \
\mu,\nu
\in M.\eqno(1.2)$$ By this definition, it can be proved directly that
for
any two \Z-modules $M,M'\subset\C$, $$\HVir\iso\Vir[M']\ \ \Rla\ \ M=a
M'
\mbox{ \ for some nonzero \ }a\in\C.\eqno(1.3)$$ Therefore we can always
suppose $1\in M$ (and in the following, in most cases, we will suppose
$1
$ is a \Z-basis element of $M$). From this, \HVir\, can be realized as
the universal central extension of the complex Lie algebra \HVB\,
consisting of the linear differential operators
$\{L_\mu=t^{\mu+1}{d\over
dt}\ |\ \mu\in M\}$ acting on the Laurant polynomial ring $L[M]=${\sl
Span}$_{\sC}\{t^\mu\ |\ \mu\in
M\}=\C[t^{b_1},...,t^{b_n},(t^{b_1})^{-1},.
 ..,(t^{b_n})^{-1}]$ if $\{b_1,...,b_n\}$ is a \Z-basis of $M$ (the
algebra \HVB\, is also called a high rank {\sl centerless} Virasoro
algebra). As the algebra $L[M]$ is actually isomorphic to the Laurant
polynomial ring $L[n]=\C[x_1,...,x_n,x^{-1}_1,...,x^{-1}_n]$ and its
derivation algebra is the Witt algebra $W[n]$, \HVB\, can also be
realized as a subalgebra of $W[n]$:
$$\HVB=\left\{p(x_1,...,x_n)\sum_{i=1}
 ^n b_ix_i{\partial\over\partial x_i}\ |\ p(x_1,...,x_n)\in
L[n]\right\}.
\eqno(1.4)$$ Here if we choose $b_1=1$, then ${d\over
dt}=x_1^{-1}\sum_{i=
1}^n b_ix_i{\partial\over\partial x_i}$. In general, if $\sum_{i=1}^nm_i
b_i=1$, for some $m_i\in\Z$, then ${d\over dt}=\prod_{i=1}^nx_i^{-m_i}
\sum_{i=1}^nb_ix_i{\partial\over\partial x_i}$. From this, we see that $
\HVB={\sl der\,}L[M]\Rla n=1$.\par Like the Virasoro algebra, \HVir\, is
$
M$-graded $$\matrix{\HVir=\OPLUS{\ }{\mu\in M}\HVir_\mu,\ \ \ \
\HVir_\mu=
\left\{^{^{\dis\C L\mu,}}_{_{\dis\C L_0\oplus\C c,}}{\ \
}^{^{\dis\mu\not=
0,}}_{_{\dis\mu=0,}}\right.\mbox{\rb{-20pt}{\,}}\hfill\cr [\HVir_\mu,
\HVir_\nu]\subset\HVir_{\mu+\nu},\ \ \mu,\nu\in M.\hfill\cr}\eqno(1.5)$$
It is worth mentioning that unlike the Virasoro algebra, when $n\ge2$,
\HVir\,is not \Z-graded unless we allow the homogeneous spaces to be
infinite dimensional.\par As in general, $M\not\subset\C$, the notion of
highest or lowest weight modules over \HVir\, will have no meaning
unless
a proper ordering on $M$ is defined. For a different ordering, a highest
or lowest weight module could look very different, and the weight spaces
of a highest or lowest weight module over \HVir\, are, in general,
infinite dimensional. Therefore, it is a hard task to study the
representation theory of the high rank Virasoro algebras.\par\ni {\bf
Definition 1.1}. (1) For any \Z-basis $B=\{b_1,...,b_n\}$ of $M$, define
\par\ni\hs{45pt} ${\sl deg}_B\mu=\sum_{i=1}^nm_i,\mbox{ \ for any \
}\mu=
\sum_{i=1}^nm_ib_i\in M,$\hfill(1.6a) \nl\hs{45pt} $M=\OPLUS{\
}{i\in\sZ}
M_B^{[i]},\ \ \ \ M_B^{[i]}=\{\mu\in M\ |\ {\sl deg}_B\mu=i\},$\hfill
(1.6b) \nl\hs{45pt} $\HVir=\OPLUS{\ }{i\in\sZ}\HVir_B^{[i]},\ \ \ \
\HVir_
B^{[i]}=\OPLUS{\ \ \ }{{\sl deg}_B\mu=i}\HVir_\mu.$\hfill(1.6c) \nl\hs
{45pt} $\HVir_{\ssc B}^{\ssc[\ge k]}={\sl Span}_{\sC}\{L_\mu\ |\
\mu=\sum_
{i=1}^nm_ib_i,m_i\ge k,i=1,...,n\},\ \ k\in\Z.$\hfill(1.6d) \par\ni(2) A
Harish-Chandra module over \HVir\, is a module $V$ such that\par\ni\hs
{45pt} $V=\OPLUS{\ }{\l\in\sC}V_\l,\ \ \ V_\l=\{v\in V\ |\ L_0v=\l v\},\
\ \ {\sl dim\,}V_\l<\INF,\ \ \ \l\in\C.\hfill(1.6e)$\nl (3) A
generalized
highest weight module over \HVir\, is an indecomposable Harish-Chandra
module generated by a generalized highest weight vector $v_\L$ such that
$
\HVir_{\ssc B}^{\ssc[\ge0]}v_\L\!=\!0$ for some \Z-basis $B$ of
$M_{_{\,}}
$. \nl(4) A uniformly bounded module over \HVir\, is a module $V$ with
weight spaces uniformly bounded, i.e., there exists some integer $N$
such
that ${\sl dim\,}V_\l\le N$ for all $\l\in\C$. \nl(5) A module of the
intermediate series over \HVir\, is an indecomposable module $V$ such
that ${\sl dim\,}V_\l\le1$ for all $\l\in\C$.\qed\par For the Virasoro
algebra, Kaplansky {\sl et al}\,\ [\ul7] classified modules of the
intermediate series over \Vir. Then Kac[\ul6] conjectured that a simple
Harish-Chandra module over \Vir\, is either a highest or lowest weight
module (note that a lowest weight module is a generalized highest weight
module by our definition), or else a module of the intermediate
series.\par This conjecture was partially proved by Chari {\sl et al}\,\
[\ul1], Martin {\sl et al}\,\ [\ul9], and the author [\ul{13}] and also
fully proved by Mathieu [\ul{10}]. \,In [\ul{14}], we were able to
extend
this result to the super-Virasoro algebras. \,In [\ul{15}], we
constructed and classified all modules of the intermediate series over
the high rank Virasoro algebras and over the high rank super-Virasoro
algebras. We proved that a module of the intermediate series over
\HVir\,
is $A_{a,b},A(a'),B(a')$ or one of their quotient modules for some $a,b
\in\C,a'\in\C\cup\{\INF\}$, where $A_{a,b},A(a'),B(a')$ all have a basis
$
\{x_\mu\ |\ \mu\in M\}$ such that $c$ acts trivially (i.e., $c x_\mu=0$
for $\mu\in M$), and $$\matrix{\hfill A_{a,b}:&L_\mu x_\nu=(a+\nu+\mu
b)x_
{\mu+\nu},\VS\hfill\cr A(a'):&L_\mu x_\nu=(\nu+\mu)x_{\mu+\nu},\nu\ne0,\
L_\mu x_0=\mu(1+(\mu+1)a')x_\mu,\VS\hfill\cr B(a'):&L_\mu x_\nu=\nu x_{
\mu+\nu},\nu\ne-\mu,\ L_\mu x_{-\mu}=-\mu(1+(\mu+1)a')x_0,\hfill\cr }
\eqno(1.7)$$ where $\mu,\nu\in M$, and here we use the same convention
as
in [\ul{14},\ul{15}] that if $a'=\INF,$ then $1+(\mu+1)a'$ is regarded
as
$\mu+1$. Note that we have $$\mbox{(i) \ }A_{a,b}\mbox{ \ is simple \ }
\Rla\ a\not\in M\mbox{ \ or \ }a\in M,b\ne0,1,\mbox{ \ and (ii) \
}A_{a,1}
=A_{a,0}\mbox{ \ if \ }a\not\in M.\eqno(1.8)$$\par It is the aim of the
present paper to extend the Kac's conjecture to the high rank Virasoro
algebras. It is interesting to see that although the modules over
\HVir\,
are in general far more complicated than those over \Vir, nevertheless,
the simple Harish-Chandra modules over \HVir\, look very similar to
those
over \Vir, i.e., we have the following main result of the present
paper.\par\ni{\bf Theorem 1.2}. A simple Harish-Chandra module over
\HVir
\, is either a generalized highest weight module, or else a module of
the
intermediate series. \par\ni{\bf Remark 1.3}. It is worth mentioning
that
although we have Theorem 1.2 for \HVir, the representation of the
generalized highest weight modules over \HVir\, is still remaining an
open problem. It may be possible and interesting to consider some
special
cases of generalized highest weight modules over \HVir.\qed\par Finally,
we point out that we are able to employ this result to classify
indecomposable Harish-Chandra modules over \HVir\, in [\ul{16}]. \par\
\par\ \par\cl{\bf\S2. Uniformly bounded modules}\par Suppose $V$ is an
indecomposable Harish-Chandra module over \HVir. For $a\in\C$, let
$V(a)=
\sum_{\mu\in M}V_{\mu+a}.$ Clearly, $V(a)$ is a submodule of \HVir, and
$
V$ can be decomposed as a direct sum of different $V(a),a\in\C/M$,
therefore we have $$V=V(a)=\sum_{\mu\in M}V_{\mu+a}\mbox{ \ \ for some \
}
a\in\C.\eqno(2.1)$$ In the rest of paper, in most cases, we will suppose
$
n\ge2$, as for $n=1$, \HVir\, is simply the (rank one) Virasoro algebra
and all the results follow from [\ul{10}]. First we have\par\ni{\bf
Lemma
2.1}. Let $V$ be a simple Harish-Chandra module over \HVir. Let
$B=\{b_1,.
 ..,b_n\}$ be as in Definition 1.1. Suppose $\HVir_{\ssc B}^{\ssc[\ge
k]}v=
0$ for some nonzero $v\in V$ and some non-negative integer $k$. Then $V$
is a generalized highest weight module.\par\ni{\bf Proof}. Let $$b'_i=
\sum_{j=1}^i(k+i-j+1)b_j+k\sum_{j=i+1}^nb_j,i=1,...,n.\eqno(2.2)$$ One
can easily check that the determinant of coefficients of $\{b_i\ |\
i=1,..
,n\}$ is 1, hence, $B'=\{b'_1,...,b'_n\}$ is also a \Z-basis of $M$.
From
(2.2), we have $\HVir_{\ssc B'}^{\ssc[\ge0]}\subset \HVir_{\ssc
B}^{\ssc[
\ge k]}$, and so $\HVir_{\ssc B'}^{\ssc[\ge0]}v=0$. Therefore $v$ is a
generalized highest weight vector and $V$ is a generalized highest
weight
module. \qed\par\ni{\bf Lemma 2.2}. A simple Harish-Chandra module $V$
over \HVir\, without a generalized highest weight vector must be
uniformly bounded. \par\ni To prove this result, we need the following
result. \par\ni{\bf Lemma 2.3}. For any $\mu=\sum_{i=1}^nm_ib_i\in M$,
let $\LL$ be the Lie subalgebra of \HVir\, generated by $A=\{L_\nu\ |\
-1
\le{\sl deg}_B(\nu-\mu)\le1\}$, then there exists a \Z-basis $B'$ of $M$
and a non-negative integer $K$ such that $\LL\supset\HVir_{\ssc
B'}^{\ssc[
\ge K]}$. \par\ni{\bf Proof}.We want to prove that there exists a \Z
-basis $B'$ of $M$ such that $\LL\supset\HVir_{\ssc B'}^{\ssc[\ge0]}$.
By
replacing $b_i$ by $-b_i$ if necessary (this does not change the set
$A$,
neither do $\LL$), we can suppose $m_i\ge0,i=1,...,n$. \nl(1) If
$m_1\not=
0,m_2\not=0$. Take
$$\matrix{b'_1=m_2\mu+b_1=(\mu+b_1)+\sum_{i=1}^{m_2-1}
\mu=(m_1m_2+1)b_1+m^2_2b_2+\sum_{i=3}^nm_2m_ib_i,\VS\hfill\cr
b'_2=m_1\mu-
b_2=(\mu-b_2)+\sum_{i=1}^{m_1-1}\mu=m^2_1b_1+(m_1m_2-1)b_2+\sum_{i=3}^nm_
1m_ib_i,\VS\hfill\cr b'_i=b'_1+b_i=(\mu+b_1+b_i)+\sum_{i=1}^{m_2-1}\mu,\
\ i=3,...,n.\hfill \cr}\eqno(2.3)$$ From the expressions (2.3), we can
easily solve $b_i,i=1,...,n,$ as an integral linear combination of $B'=
\{b'_1,...,b'_n\}$. Thus $B'$ is a \Z-basis of $M$. We also see that
$\{L_
{b'_i}\ |\ i=1,...,n\}\subset \LL$. For example, from
$$[L_\mu,...,[L_\mu,
L_{\mu+b_1}]...]\ (m_2-1 \mbox{ copies of }L_\mu)=aL_{b'_1},\mbox{ \
where \ }a=\prod_{i=0}^{m_2-2}(i\mu+b_1)\ne0,\eqno(2.4)$$ and the {\sl
l.h.s.} of (2.4) is in $\LL$, we see $L_{b'_1}\in\LL$. This proves that
$
\HVir_{\ssc B'}^{\ssc[\ge0]}$, which is the Lie subalgebra generated by
$
\{L_{b'_1}\ |\ i=1,...,n\}$, is contained in $\LL$.\nl(2) If $m_i=0$ for
some $i,1\le i\le n$. Say, $m_1=0$. Take
$$b'_1=\mu+b_1,b'_i=b'_1+b_i=\mu+
b_1+b_i,\ i=2,...,n.$$ Then as above $B'=\{b'_1,...,b'_n\}$ is a
\Z-basis
of $M$ and $\HVir_{\ssc B'}^{\ssc[\ge0]}\subset\LL$. \qed\par\ni{\bf
Proof of Lemma 2.2}. Let $V$ be as in (2.1). For any $\mu\in M$, let $A$
and $\LL$ be as in Lemma 2.3. By Lemma 2.1, we see that for any $0\not=v
\in V$, we have $\LL v\not=0$, thus, $Av\not=0$. This means that
$\cap_{-
1\le{\sl deg}_B(\nu-\mu)\le1}{\sl Ker}(L_\nu|_V)=0$. In particular, $
\oplus_{-1\le{\sl deg}_B(\nu-\mu)\le1}L_\nu|_{V_{a-\mu}}:V_{a-\mu}\rar
\oplus_{-1\le{\sl deg}_B\nu\le1}V_{a+\nu}$ is an injection. Therefore,
${
\sl dim\,}V_{a-\mu}\le N$, where $N=\sum_{-1\le{\sl deg}_B\nu\le1}{\sl
dim\,}V_{a+\nu}$ is a fixed integer. As $\mu\in M$ is arbitrary, this
proves that $V$ is uniformly bounded.\qed\par \ni{\bf Lemma 2.4}. Let
$V$
as in (2.1) be a uniformly bounded module over \HVir. There exists a
non-negative integer $N$ such that ${\sl dim\,}V_{a+\mu}=N$ for all $\mu
\in M,\,a+\mu\ne0$. \par\ni{\bf Proof}. By decomposing $V$ as a direct
sum of indecomposable submodules, we see that it suffices to prove the
result when $V$ is indecomposable.\par If $n=1$, the result follows from
[\ul9]. Suppose now $n\ge2$. Let $B=\{b_1,...,b_n\}$ be a \Z-basis of $M
$. For $i=1,...,n$, let \Vir$_i$ be the Virasoro algebra generated by $
\{L_{kb_i}\ |\ k\in\Z\}$. For any $\mu\in M$, let
$V(\mu,i)=\sum_{k\in\Z}
V_{a+\mu+kb_i}$. Clearly $V(\mu,i)$ is a uniformly bounded module over
\Vir$_i$, therefore we have ${\sl dim\,}V_{a+\mu}={\sl
dim\,}V_{a+\mu+kb_
i}$, for all $\mu\in M,k\in\Z,i=1,...,n$, such that
$a+\mu\ne0,a+\mu+kb_i
\ne0.$ From this, we can conclude that ${\sl dim\,}V_{a+\mu}$ is a
constant for all $\mu\in M,\,a+\mu\ne0$.\qed\par\ \par\ \par \cl{\bf\S3.
Proof of Theorem 1.2}\par Now we are in a position to prove Theorem 1.2.
Suppose $V$ is a simple Harish-Chandra module over \HVir. Then $V$ has
the form (2.1) and by Lemma 2.2 and Lemma 2.4, we can suppose $V$ is
uniformly bounded such that ${\sl dim\,}V_{a+\mu}=N$ for all $\mu\in M,
\,a+\mu\ne0$.\par If $N=0$, $V$ is just a trivial module, and the result
is obvious true. Therefore, we suppose $N\ge1$. If $n=1$, the result
follows from [\ul{10}]. Suppose $n\ge 2$, and by inductive assumption,
we
can suppose the result holds for any rank $n-1$ Virasoro algebra. By the
statement after (1.3), we can take a \Z-basis
$B=\{1=b_1,b_2,...,b_{n-1},
b_n=d\}$ of $M$, so $$M=M_1\oplus\Z d,\ \
M_1=\Z\oplus(\oplus_{i=2}^{n-1}
\Z b_i).\eqno(3.1)$$ Let \HVira\, be the rank $n-1$ Virasoro algebra
generated by $\{L_\mu\ |\ \mu\in M_1\}$. For $k\in\Z$, let
$V^{(k)}=\sum_
{\mu\in M_1}V_{a+kd+\mu}$, then it is a uniformly bounded \HVira
-module.\par As $B$ is a \Z-basis of $M$, one sees that there exists at
most one $k\in\Z$ such that $a+kd\in M_1$. Fix a $k_0$ such that $$a+kd
\not\in M_1\mbox{ \ for \ }k\ge k_0.\eqno(3.2)$$ Now first we take
$k=k_0
$. By inductive assumption, the result of Theorem 1.2 holds for the rank
$
n-1$ Virasoro algebra \HVira, thus there exists a simple submodule
$V^{(k,
1)}$ of $V^{(k)}$ over \HVira. As $V$ is simple, the center element $c$
acting on $V$ must be a scalar, by induction on $n$, we can suppose $c$
acts trivially on $V^{(k,1)}$, thus $c$ must also trivially act on $V$.
By (1.8), $V^{(k,1)}$ has the form $A_{a+kd,b}$, i.e., there exists a
basis $\{x_\mu\ |\ \mu\in M_1\}$ of $V^{(k,1)}$, such that $$L_\mu
x_\nu=(
\ol\nu+\mu b)x_{\mu+\nu},\mu,\nu\in M_1,\eqno(3.3)$$ where, here and the
following, to simplify the notation, we take the following convention,
$$
\ol\nu=a+kd+\nu,\mbox{ \ for \ }\nu\in M_1 \mbox{ (when $k$ is fixed)}.
\eqno(3.4)$$ Furthermore, we can suppose that among all simple modules
of
$V^{(k)}$, $V^{(k,1)}$ is chosen to be one such that ${\sl Re}(b)$ is
minimum, where ${\sl Re}(b)$ is the real part of the complex number $b$.
We can take a composition series of $V^{(k+1)}$ over \HVira:
$$0=V^{(k+1,
0)}\subset V^{(k+1,1)}\subset...\subset V^{(k+1,N)}=V^{(k+1)},\eqno(3.5)
$$ such that for $i=1,...,N$, $$V^{(k+1,i)}/V^{(k+1,i-1)}\mbox{ \ has
the
form \ }A_{a+(k+1)d,b'}\mbox{ \ for some \ }b'\in\C,\eqno(3.6)$$ where $
b'$ may depend on $i$. Clearly we have $$L_d V^{(k,1)}\subset L_d
V^{(k)}
\subset V^{(k+1)}.\eqno(3.7)$$ {\bf Lemma 3.1}. $L_d
V^{(k,1)}\not=0.$\par
\ni{\bf Proof}. Suppose $L_d V^{(k,1)}=0.$ Note that for any $\mu\in
M_1,
m\in\Z_+\bs\{0\}$, $L_{\mu+md}$ can be generated by $\HVira\cup\{L_d\}$.
By this and (3.3), we obtain by induction on $m$ that
$L_{\mu+md}V^{(k,1)}
=0$ if $m>0$. Let $W$ be the \HVir-submodule of $V$ generated by
$V^{(k,1)
}$. Since, as space, $\HVir=\HVir_1\oplus\HVir_2$, where $\HVir_1$ and $
\HVir_2$ are the Lie subalgebra of \HVir\, generated by
$\HVira\cup\{L_{-
d}\}$ and respectively by $\{L_{\mu+md}\ |\ \mu\in
M_1,m\in\Z_+\bs\{0\}\,
\}$, by decomposing the universal enveloping algebra of \HVir\, as $U(
\HVir)=U(\HVir_1)U(\HVir_2)$, we have
$W=U(\HVir_1)U(\HVir_2)V^{(k,1)}=U(
\HVir_1)V^{(k,1)}$. From this, we obtain that $W_{\ol\mu+md}=0$ if $\mu
\in M_1,m>0$. On the other hand, $W_{\ol\mu}\supset
V^{(k,1)}_{\ol\mu}=\C
x_\mu,\mu\in M_1$. Thus Lemma 2.4 is violated by the \HVir-module $W$.
This contradiction proves that $L_d V^{(k,1)}\not=0.$\qed\par By (3.7)
and (3.5), let $N_1\ge0$ be the integer such that $$L_d V^{(k,1)}\not
\subset U,\mbox{ \ but \ }L_d V^{(k,1)}\subset U',\mbox{ \ where \
}U=V^{(
k+1,N_1)},U'=V^{(k+1,N_1+1)}.\eqno(3.8)$$ Furthermore, we can choose
(3.5) to be a composition series such that $N_1$ is minimum, among all
composition series of $V^{(k+1)}$. By (3.6), we can take a basis
$\{y_\mu
\ |\ \mu\in M_1\}$ of $U'/U$ such that $$\matrix{L_\mu
y_\nu\equiv(\ol\nu+
d+\mu b')y_{\mu+\nu}\,({\rm mod\,}U),\ \mu,\nu\in M_1.\VS\hfill\cr L_d
x_
\nu\equiv a_\nu y_\nu\,({\rm mod\,}U),\mbox{ \ for some \ }a_\nu\in\C,\
\nu\in M_1.\hfill\cr }\eqno(3.9)$$ {\bf Lemma 3.2}. We have $b'=b$ or
$b-
1$. Furthermore by rescaling $y_\nu,\nu\in M_1$ if necessary, we have
$$a_
\nu=\left\{\matrix{\ol\nu+bd&\mbox{if \ }b'=b,\VS\hfill\cr 1&\mbox{if \
}
b'=b-1.\hfill\cr }\right.\eqno(3.10)$$ \par\ni{\bf Proof}. For $\mu,\mu'
\in M_1$, by $[L_\mu,[L_{\mu'},L_d]]=(d-\mu')(d+\mu'-\mu)L_{\mu+\mu'+d}$
and $[L_{\mu+\mu'},L_d]=(d-\mu-\mu')L_{\mu+\mu'+d}$, we have $$\matrix{
\hfill(d-\mu-\mu')(L_\mu L_{\mu'}L_d-L_\mu L_d L_{\mu'}-L_{\mu'}L_d
L_\mu+
L_d L_{\mu'}L_\mu)&\VS\cr\hfill -(d-\mu')(d+\mu'-\mu)(L_{\mu+\mu'}L_d-L_
dL_{\mu+\mu'})&\!\!\!\!=0.\cr }\eqno(3.11)$$ Apply (3.11) to $x_\nu$, by
(3.9) and (3.3), we obtain (where the calculation is conducted under ${
\rm mod\,}U$ and we count the coefficients of $y_\nu,\nu\in M_1$)\par\ni
\hs{5pt}$\matrix{(d-\mu-\mu'){\bf(}(\ol\nu+d+\mu'b')(\ol\nu+\mu'+d+\mu
b')
a_\nu -(\ol\nu+\mu'b)(\ol\nu+\mu'+d+\mu b')a_{\nu+\mu'}\hfill&\VS\cr
\hfill-(\ol\nu+\mu b)(\ol\nu+\mu+d+\mu'b')a_{\nu+\mu}+(\ol\nu+\mu b)(\ol
\nu+\mu+\mu'b)a_{\nu+\mu+\mu'}{\bf)}&\VS\cr \hfill-(d-\mu')(d+\mu'-\mu){
\bf(}(\ol\nu+d+(\mu+\mu')b')a_\nu-(\ol\nu+(\mu+\mu')b)a_{\nu+\mu+\mu'}{
\bf)}&\!\!\!\!=0.\cr }\hfil(3.12)$\par\ni Now in (3.12), (i) replacing $
\mu'$ by $\mu$, and $\nu$ by $\nu-\mu$; (ii) replacing $\mu'$ by $-\mu$;
(iii) replacing both $\mu$ and $\mu'$ by $-\mu$, and $\nu$ by $\nu+\mu$,
we obtain three equations concerning $a_{\nu-\mu},a_\nu,a_{\nu+\mu}$
correspondingly:\par\ni\hs{5pt}$\matrix{\hfill((d-2\mu)(\ol\nu-\mu+d+\mu
b')(\ol\nu+d+\mu b')-d(d-\mu)(\ol\nu-\mu+d+2\mu b'))a_{\nu-\mu}&\VS\cr
\hfill-2(d-2\mu)(\ol\nu-\mu+\mu b)(\ol\nu+d+\mu b')a_\nu&\VS\cr
\hfill+((
d-2\mu)(\ol\nu-\mu+\mu b)(\ol\nu+\mu b)+d(d-\mu)(\ol\nu-\mu+2\mu b))a_{
\nu+\mu}&\!\!\!\!=0,\cr }\hfill\!\!\!(3.13a)$\par\ni\hs{5pt}$\matrix{
\hfill(\ol\nu-\mu b)(\ol\nu-\mu+d+\mu b')a_{\nu-\mu}&\VS\cr \hfill-((\ol
\nu+d-\mu b')(\ol\nu-\mu+d+\mu b')+(\ol\nu+\mu b)(\ol\nu+\mu-\mu b)-(d+
\mu)(d-2\mu))a_\nu&\VS\cr \hfill+(\ol\nu+\mu b)(\ol\nu+\mu+d-\mu b')a_{
\nu+\mu}&\!\!\!\!=0,\cr }\hfill\!\!\!(3.13b)$\par\ni\hs{5pt}$\matrix{
\hfill((d+2\mu)(\ol\nu+\mu-\mu b)(\ol\nu-\mu b)+d(d+\mu)(\ol\nu+\mu-2\mu
b))a_{\nu-\mu}&\VS\cr\hfill -2(d+2\mu)(\ol\nu+\mu-\mu b)(\ol\nu+d-\mu
b')
a_\nu &\VS\cr\hfill((d+2\mu)(\ol\nu+\mu+d-\mu b')(\ol\nu+d-\mu b')-d(d+
\mu)(\ol\nu+\mu+d-2\mu
b'))a_{\nu+\mu}&\!\!\!\!=0.\cr}\hfill\!\!\!(3.13c)
$\par\ni First suppose that the set $\{\nu\in M_1\ |\ a_\nu\ne0\}$ is
finite. Then for any $\nu\in M_1$, we can always choose $\mu\in M_1$
such
that $a_{\nu\pm\mu}=0$ and such that the coefficient of $a_\nu$ in
(3.13a) is nonzero. This shows that $a_\nu=0$ for all $\nu\in M_1$, and
so (3.9) tells us that $L_d V^{(k,1)}\subset U$. This contradicts (3.8).
Therefore, we must have that $$\mbox{The set \ }\{\nu\in M_1\ |\ a_\nu
\ne0\}\mbox{ \ is infinite.}\eqno(3.14)$$ This means that for any
$\nu\in
M_1$, there exist an infinite number of $\mu\in M_1$ such that the
determinant of coefficients of $a_{\nu-\mu},a_\nu,a_{\nu+\mu}$ in
(3.13a-c) is zero. Denote this determinant by $\D(\ol\nu,\mu)$. As $\D(
\ol\nu,\mu)$ is a polynomial on $\mu$ (when $\nu$ is fixed), this shows
that for any $\nu\in M_1$, $D(\ol\nu,\mu)=0$ for all $\mu\in\C$.
Similarly as $D(\ol\nu,\mu)$ is also a polynomial on $\ol\nu$ (when
$\mu$
is fixed). Thus, we have $$D(\ol\nu,\mu)=0\mbox{ \ for all \ }\nu,\mu\in
\C.\eqno(3.15)$$ By a little lengthy calculation (or simply using a
computer), we obtain
$$\matrix{D(\ol\nu,\mu)=&\!\!\!\!(b'-b)(1+b'-b)\mu^6
\times\VS\hfill\cr&\matrix{(4d(3(b+b')^2-7b-5b'+2)\ol\nu\VS\hfill\cr
\,\,
\,+4(b+b')((b+b')^2(b-b')-2(b+b')(b-2b')-(b+5b')+2)\mu^2\VS\hfill
\cr\,\,
\,\matrix{+d^2(-(b+b')^3(b-b')-2(b+b')(b^2+3bb'+4b'^2)\VS\,\,\,\cr
\hfill+
19b^2+34bb'+19b'^2-24b-16b'+4))\cr}\hfill\cr}\hfill\cr }\eqno(3.16)$$ By
(3.15), (3.16), we can solve that $$b'=b,\mbox{ or }b'=b-1,\mbox{ or
}b=1,
b'=-1,\mbox{ or }b=0,b'=1,\mbox{ or }b=2,b'=0.\eqno(3.17)$$ (1) If
$b'=b$
or $b'=b-1$. It is straightforward to verify that (3.10) is a solution
to
(3.13a-c). Note that the determinant of the coefficients of
$a_{\nu-\mu}$
and $a_\mu$ in (3.13a\&b) is non-zero. This means that we can solve $a_{
\nu-\mu}$ and $a_\mu$ in terms of $a_{\nu+\mu}$, thus the solution to
(3.13a-c) is unique up to scalars. By rescaling $y_\nu,\nu\in M_1$ if
necessary, we can suppose that the solution to (3.13a-c) is unique.\nl
(2) If $b=1,b'=-1$. By (1.8), $A_{a+kd,1}\iso A_{a+kd,0}$, by our choice
of $b$ such that ${\sl Re}(b)$ is minimum, in this case, we should have
$
b=0$. Thus this becomes a special case of (1).\nl (3) Suppose $b=0,b'=1$
or $b=2,b'=0$. Again from $A_{a+kd,1}\iso A_{a+kd,0}$, we can choose
$b'$
to be 0 or 1, so that this becomes a special case of (1). This completes
the proof of the lemma.\qed\par \ni{\bf Lemma 3.3}. If $b'=b-1$, then $
\{L_d x_\nu\ |\ \nu\in M_1\}$ generates a simple submodule $W$ of
$V^{(k+
1)}$. \par\ni{\bf Proof.} We need to prove in (3.8), $N_1=0$, i.e., $U=0
$. Suppose $U\not=0$. Take a composition series of $U$: $$0=U_0\subset
U_
1\subset...\subset U_{N_1}=U.\eqno(3.18)$$ Taking a basis $\{z_\nu\ |\
\nu\in M_1\}$ of $U/U_{N_1-1}$, and replacing $y_\nu$ by $L_d x_\nu$, $
\nu\in M_1$, by (3.3), (3.9) and (3.10), we now have $$\matrix{L_\mu x_
\nu=(\ol\nu+\mu b)x_{\mu+\nu},\hfill&L_d x_\nu=y_\nu,\VS\hfill\cr L_\mu
y_
\nu\equiv(\ol\nu+d+\mu(b-1))y_{\mu+\nu}+c_{\mu,\nu}z_{\mu+\nu},\hfill&L_
\mu z_\nu\equiv(\ol\nu+d+\mu b'')z_{\mu+\nu},\hfill\cr }\ \ \mu,\nu\in
M_
1,\eqno(3.19)$$ for some $b'',c_{\mu,\nu}\in\C$ (such that $c_{0,\nu}=0
$), where, here and below, the calculation is conducted under ${\sl
Mod\,}
U_{N_1-1}$. If $c_{\mu,\nu}=0$ for all $\mu,\nu\in M_1$, then (3.18) and
(3.19) tell us that we can choose a composition series (3.5) of
$V^{(k+1)}
$ with $U$ replaced by $U_{N_1-1}$, this, contradicts our choice of
$N_1$
being minimum. Thus $c_{\mu,\nu}\not=0$ for some $\mu,\nu\in M_1$. If
necessary, by a suitable choice of basis of $M_1$ in (3.1), we can
suppose $c_{jd_i,\nu}\not=0$ for some $j,i,\nu$. Without loss of
generality, we can suppose $i=1$, i.e., $d_i=1$, thus $$c_{j,\mu}\not=0
\mbox{ \ for some \ }j\in\Z,\mu\in M_1.\eqno(3.20)$$ From (3.19), we
obtain (we simply denote $c_{1,\nu}$ by $c_\nu$) $$\matrix{\hfill
L_{1+d}
x_\nu=&\!\!\!\!{1\over d-1}[L_1,L_d]x_\nu\equiv y_{\nu+1}+{c_\nu\over
d-1}
z_{\nu+1},\mbox{\rb{-15pt}{\,}}\hfill\cr \hfill y_\nu=&\!\!\!\!L_d
x_\nu=
{1\over d+2}[L_{-1},L_{1+d}]x_\nu\equiv \VS\hfill\cr&\!\!\!\!y_\nu+{1
\over d+2}(c_{-1,\nu+1}-{1\over
d-1}((\ol\nu-b)c_{\nu-1}-(\ol\nu+1+d-b'')
c_\nu))z_\nu,\mbox{\rb{-15pt}{\,}}\hfill\cr \hfill
2(\ol\nu+d)y_\nu=&\!\!
\!\![L_{-1},L_1]y_\nu\equiv\VS\hfill\cr &\!\!\!\!2(\ol\nu+d)y_\nu+((\ol
\nu+d+b-1)c_{-1,\nu+1}+(\ol\nu+1+d-b'')c_\nu \VS\hfill\cr
&\hfill-(\ol\nu+
d-b+1)c_{\nu-1}-(\ol\nu-1+d+b'')c_{-1,\nu})z_\nu.\cr }\eqno(3.21)$$ By
the second and the third equations, we obtain, for $\nu\in M_1$, $$
\matrix{c_{-1,\nu}={1\over d-1}((\ol\nu-1-b)c_{\nu-2}-(\ol\nu+d-b'')c_{
\nu-1}),&\mbox{\rb{-15pt}{\,}}\hfill\cr
(\ol\nu+b)(\ol\nu+1+d-b'')c_\nu-(
2\ol\nu^2+(2d-1)\ol\nu-d+b''-b''^2 -b^2+1)c_{\nu-1}&\VS\hfill\cr
\hfill+(
\ol\nu+d-1+b'')(\ol\nu-1-b)c_{\nu-2}&\!\!\!\!=0.\hfill\cr }\eqno(3.22)$$
Similarly, we have \par\ni\hs{3pt}$\matrix{L_{-1+d}x_\nu={1\over
d+1}[L_{-
1},L_d]x_\nu\equiv y_{\nu-1}+{1\over(d+1)(d-1)}((\ol\nu-1-b)c_{\nu-2}-(
\ol\nu+d-b'')c_{\nu-1})z_{\nu-1},\mbox{\rb{-15pt}{\,}}\hfill\cr
L_{2+d}x_
\nu={1\over d-2}[L_2,L_d]x_\nu\equiv y_{\nu+2}+{1\over d-2}c_{2,\nu}z_{
\nu+2},\mbox{\rb{-15pt}{\,}}\hfill\cr L_{2+d}x_\nu={1\over
d}[L_1,L_{1+d}]
x_\nu\equiv y_{\nu+2}+{1\over d(d-1)}((k+1+d+b'')c_\nu-(\ol\nu+1+b-d)c_{
\nu+1})z_{\nu+2},\mbox{\rb{-15pt}{\,}}\hfill\cr L_{1+d}x_\nu={1\over
d-3}[
L_2,L_{-1+d}]x_\nu\equiv y_{\nu+1}+{1\over
d-3}(c_{2,\nu-1}+{1\over(d+1)(
d-1)}(((\ol\nu-1-b)c_{\nu-2}\VS\hfill\cr \hfill-(\ol\nu\!\!+\!\!d\!\!-\!
\!b'')c_{\nu-1})(\ol\nu\!-\!1\!+\!d\!+\!2b'')\!-\!(\ol\nu\!+\!2b)((\ol\nu
\!+\!1\!-\!b)c_\nu\!-\!(\ol\nu\!+\!2\!+\!d\!-\!b'')c_{\nu+1})))z_{\nu+1}.
\cr}\!\!\!\!\hfill(3.23)$\par \ni From this, we obtain\par\ni\hs{3pt}$
\matrix{c_{2,\nu}={d-2\over
d(d-1)}((\ol\nu+1+d+b'')c_\nu-(\ol\nu+b-d+1)c_
{\nu+1}),\mbox{\rb{-15pt}{\,}}\hfill\cr d(\ol\nu\!-\!1\!+\!2b)(\ol\nu\!+
\!1\!+\!d\!-\!b'')c_\nu
\!-\!(d\ol\nu^2\!+\!(d^2\!+\!(b\!-\!2)d\!-\!2)\ol
\nu \!+\!(b\!-\!2)d^2\!-\!2b^2d\!+\!2\!-\!2b)c_{\nu-1}\VS\hfill\cr\hfill
-
(d\ol\nu^2+(d^2+(b''-2)d+2)\ol\nu-d^2-(2b''^2-b''-3)d+2b''-2)c_{\nu-2}\VS
\cr\hfill
+d(\ol\nu-2-b)(\ol\nu-2+d+2b'')c_{\nu-3}=0.\!\!\!\!\!\!\!\!\!\!
\cr }\hfill(3.24)$\par\ni Now from (3.22) and (3.24), we can eliminate
$c_
{\nu-i},i=1,2,3$, as below to obtain $p(\ol\nu)c_\nu=0$, where
$p(\ol\nu)
$ is a polynomial on $\ol\nu$: if we simply write the second equations
of
(3.22) and (3.24) as $\sum_{i=0}^2 s_i(\ol\nu)c_{\nu-i}=0$ and
$\sum_{i=0}
 ^3 t_i(\ol\nu)c_{\nu-i}=0$, by taking
$$\matrix{u_1(\ol\nu)=&\!\!\!\!\!d(
\ol\nu-1+2b)s_1(\ol\nu)-(\ol\nu+b)t_1(\ol\nu),\VS\hfill\cr u_2(\ol\nu)=&
\!\!\!\!\!d(\ol\nu-1+2b)s_2(\ol\nu)-(\ol\nu+b)t_2(\ol\nu),\VS\hfill\cr
u_
3(\ol\nu)=&\hfill-(\ol\nu+b)t_3(\ol\nu),\cr}\eqno(3.25)$$ we can
iliminate $c_\nu$ to obtain $\sum_{i=1}^3 u_i(\ol\nu)c_{\nu-i}=0$; by
taking $$\matrix{v_1(\ol\nu)=(\ol\nu-2-d+b'')u_1(\ol\nu)+d(\ol\nu-2+d+
2b'')(\ol\nu+b)s_0(\ol\nu-1),\VS\hfill\cr
v_2(\ol\nu)=(\ol\nu-2-d+b'')u_2(
\ol\nu)+d(\ol\nu-2+d+2b'')(\ol\nu+b)s_1(\ol\nu-1),\hfill\cr}\eqno(3.26)$$
we can iliminate $c_{\nu-3}$ to obtain
$$v_1(\ol\nu)c_{\nu-1}+v_2(\ol\nu)
c_{\nu-2}=0;\eqno(3.27)$$ by taking
$$\matrix{w_0(\ol\nu)=s_0(\ol\nu)v_2(
\ol\nu),\VS\hfill\cr w_1(\ol\nu)=s_1(\ol\nu)v_2(\ol\nu)-s_2(\ol\nu)v_1(
\ol\nu),\hfill\cr }\eqno(3.28)$$ we obtain
$$w_0(\ol\nu)c_\nu+w_1(\ol\nu)
c_{\nu-1}=0;\eqno(3.29)$$ Finally by taking
$$p(\ol\nu)=v_1(\ol\nu+1)w_1(
\ol\nu)-v_2(\ol\nu+1)w_0(\ol\nu),\eqno(3.30)$$ we obtain
$p(\ol\nu)c_\nu=
0$. Now we claim that $c_\nu=0$ for all $\nu\in M_1$. If not, from
(3.22), we see that there must be an infinite number of $\nu$ such that
$
c_\nu\not=0$, and such that $p(\ol\nu)=0$ correspondingly. Thus
$p(\ol\nu)
=0$ for all $\nu\in\C$. From above, by a lengthy computation or using a
computer, we have
$$p(\ol\nu)=-d(d+1)(b''-b)(b''-b-1)(\ol\nu+b)(\ol\nu+d-
1+b'')(p_0\ol\nu^2+p_1\ol\nu+p_2),\eqno(3.31{\rm a})$$ for some
$p_0,p_1,
p_2\in\C$, where $p_0=(b+b'')(b+b''-1)$, and $p_1$ can be written as
$$p_
1=p'_1p_0-2(3b''^2-3b''-1)bd-2(3b''^3-4b''^2-b''+1)d,\eqno(3.31{\rm
b})$$
for some $p'_1\in\C$. Thus we can solve $b''=b$ or $b''=b+1$ or $b=\pm1,
b''=\,^{_-}_{^+}1$ or $b=1,b''=0$. \nl (1) Suppose $b''=b$. By the
second
equation of (3.22), (3.27) and (3.29), we get $$c_\nu=c_{\nu-1},\nu\in
M_
1.\eqno(3.32)$$ From this, by (3.22) and (3.24), we get
$$c_{i,\nu}={i(d-
i)\over d-1}c_\nu,\eqno(3.33)$$ for $i=0,\pm1,2$. Now using
$[L_1,L_{i-1}]
=(i-2)L_i$, $[L_{-i+1},L_{-1}]=(i-2)L_{-i}$, applying them to $y_\nu$,
using induction on $i$, we obtain that (3.33) holds for all $i\in\Z$.
Applying $[L_i,L_j]=(j-i)L_{i+j}$ to $y_\nu$, we can finally obtain $c_
\nu=0$. \nl(2) Suppose $b''=b+1$. As in (1), we can obtain in this case
$$
c_{i,\nu}={i(d-i)\over d-1}(\ol\nu+(2b+1)d+i(b+1))c'_\nu,\eqno(3.34)$$
where $c'_\nu=c'_{\nu-1}$ for $\nu\in M_1$. Again, we can get
$c'_\nu=0$.
This proves that $c_\nu=0$ for all $\nu\in M_1$. \nl (3) Suppose
$b=\pm1,
b''=\,^{_-}_{^+}1$ or $b=1,b''=0$. As above, we have $$\left\{\matrix{c_
\nu={c'_\nu\over\ol\nu(\ol\nu+1)},\hfill&\mbox{and \ }c_{i,\nu}={i(d-i)
\over d-1}{c'_\nu\over\ol\nu(\ol\nu+i)},\hfill &\mbox{if \ }b=1,b''=-1,
\VS\hfill\cr c_\nu=(\ol\nu-d)(\ol\nu+d+1)c'_\nu,\hfill&\mbox{and \
}c_{i,
\nu}={i(d-i)\over d-1}(\ol\nu-d)(\ol\nu+d+i)c'_\nu,\hfill &\mbox{if \
}b=-
1,b''=1,\VS\hfill\cr c_\nu={c'_\nu\over\ol\nu+d+1},\hfill&\mbox{and \
}c_
{i,\nu}={i(d-i)\over d-1}{c'_\nu\over\ol\nu+d+i},\hfill &\mbox{if \
}b=1,
b''=0,\hfill\cr}\right.$$ where $c'_\nu=c'_{\nu-1},\ \nu\in M_1$. Again,
we can obtain $c'_\nu=0$, thus $c_\nu=0$ for all $\nu\in M_1$.\par In
all
cases, by (3.22) and (3.24) and by induction, we obtain that
$c_{i,\nu}=0
$ for all $i\in\Z,\nu\in M_1$. This contradicts (3.20). This completes
the proof of the lemma.\qed\par\ni {\bf Lemma 3.4}. If $b'=b$, then
$\{L_
d x_\nu\ |\ \nu\in M_1\}$ generates a simple submodule $W$ of $V^{(k+1)}
$. \par\ni{\bf Proof.} First suppose $\ol\nu+bd\not=0$ for all $\nu\in
M_
1$. In this case, as in the proof of Lemma 3.3, we can replace $y_\nu$
by
${1\over\ol\nu+bd}L_dx_\nu$. Then we have the first and last equations
of
(3.19) and $$L_d x_\nu=(\ol\nu+bd)y_\nu,\ L_\mu y_\nu\equiv(\ol\nu+d+\mu
b)y_{\mu+\nu}+c_{\mu,\nu}z_{\mu+\nu},\ \ \mu,\nu\in M_1.\eqno(3.35)$$
Now
the proof is exactly analogous to, though more involved than, that of
Lemma 3.3. Thus we can obtain
\par\ni\hs{2pt}$\matrix{(\ol\nu\!-\!1\!+\!(
d\!+\!1)b)c_{-1,\nu}={1\over d-1}((\ol\nu\!-\!1\!-\!b)(\ol\nu\!-\!2\!+
\!bd)c_{\nu-2}-(\ol\nu\!-\!1\!+\!bd)(\ol\nu\!+\!d\!-\!b'')c_{\nu-1},\VS
\hfill\cr (\ol\nu\!+\!bd)c_{2,\nu}\!=\!{d\!-\!2\over d(d\!-\!1)}({\sc(}
\ol\nu\!+\!bd{\sc)}{\sc(}\ol\nu\!+\!1\!+\!d\!+\!b''{\sc)}c_\nu
\!-\!({\sc(
}\ol\nu\!+\!b{\sc)}{\sc(}\ol\nu\!+\!1\!+\!bd{\sc)}\!-\!{\sc(}d\!-\!1{\sc)}
{\sc(}\ol\nu\!+\!{\sc(}d\!+\!1{\sc)}b{\sc)})c_{\nu+1}).\hfill\cr
}\hfill(
3.36)$\par\ni Similarly, we have $p(\ol\nu)c_\nu=0$ for some polynomial
$
p(\ol\nu)$. Use a computer's help, we can compute $p(\ol\nu)$ and find
out that $p(\ol\nu)\not=0$. This forces $c_\nu=0$.\par Suppose now there
exists $\nu_0\in M_1$ such that $\ol{\nu_0}+bd=0$. In this case, we must
have $b\not=1$, otherwise by (3.3), $0=\ol{\nu_0}+bd=a+(k+1)d+\nu_0$, a
contradiction with (3.2). Now we take $y_\nu={1\over\ol\nu+bd}L_d x_\nu$
if $\nu\not=\nu_0$ and take $y_{\nu_0}={1\over\ol\nu_0-1+d+b}L_1
y_{\nu_0-
1}$ (we have $\ol\nu_0-1+d+b=-1+d+b-bd=(d-1)(1-b)\not=0$ as by (3.1), $d
\not=1$), then we have (3.35) with the first equation replaced by $$L_d
x_
\nu\equiv(\ol\nu+bd)y_\nu+\d_{\nu,\nu_0}c z_\nu,\ \ \nu\in
M_1,\eqno(3.37)
$$ for some $c\in\C$. By our choice of $y_{\nu_0}$, we have
$c_{\nu_0-1}=
c_{1,\nu_0-1}=0$. Now we can proceed as above to prove that $c_\nu=0$
for
all $\nu\in M_1$ and that $c=0$.\qed\par By Lemma 3.3 and Lemma 3.4, we
can now let $V^{(k+1,1)}$ be the simple \HVira-submodule of $V^{(k+1)}$
generated by $V^{(k,1)}$. Now we claim that $b'=b$. If not, by taking
$d'=
-d$, and considering $L_{d'}V^{(k+1,1)}$. Following exactly the same
arguments as before, we then have a simple \HVira-submodule $V'^{(k,1)}$
of $V^{(k)}$, such that $L_{d'}V^{(k+1,1)}\subset V'^{(k,1)}$ and
$V'^{(k,
1)}\cong A_{a+kd,b_1}$ with $b_1=b'$ or $b_1=b'-1$. In particular, ${\sl
Re}(b_1)<{\sl Re}(b)$. This contradicts our choice of $b$ that ${\sl
Re}(
b)$ is minimum. Therefore $b'=b$. Now for $i=1,2$, let $\wt V^{(k+i,1)}$
be the sum of all simple \HVira-submodules of $V^{(k+i)}$ which are
isomorphic to $A_{a+(k+i)d,b}$, then the above arguments show that $L_d
\wt V^{(k,1)}\subset\wt V^{(k+1,1)}$ and that $L_{-d}\wt V^{(k+1,1)}
\subset\wt V^{(k,1)}$. Suppose $\wt V^{(k+i,1)}$ has $r_i$ composition
factors ($i=1,2$), then by Lemma 3.1, we must have $r_1=r_2=r$. Now
choose a suitable basis $\{x^{(i)}_\nu\ |\ \nu\in M_1,i=1,...,r\}$ of $
\wt V^{(k,1)}$ such that $x^{(1)}_0$ is an eigenvector of $L_{-d}L_d$,
and such that $$L_\mu X_\nu=(\ol\nu+\mu b)X_{\mu+\nu},\mu,\nu\in M_1,
\eqno(3.38)$$ where, $X_\nu=(x^{(1)}_\nu,...,x^{(r)}_\nu)^t$ for $\nu\in
M_1$. Lemma 3.2 allows us to choose a basis $\{y^{(i)}_\nu\ |\ \nu\in
M_1,
i=1,...,r\}$ of $\wt V^{(k+1,1)}$ such that $$\matrix{L_\mu
Y_\nu=(\ol\nu+
d+\mu b)Y_{\mu+\nu},\mu,\nu\in M_1,\VS\hfill\cr L_d X_\nu=(\ol\nu+bd)Y_
\nu,\nu\in M_1,\hfill\cr }\eqno(3.39)$$ where,
$Y_\nu=(y^{(1)}_\nu,...,y^
{(r)}_\nu)^t$ for $\nu\in M_1$. Now suppose $L_{-d}Y_\nu=A_\nu X_\nu$
for
some $r\times r$ complex matrix $A_\nu$, $\nu\in M_1$. Then as in the
proof of Lemma 3.2, we can solve that $A_\nu=(\ol\nu+d-bd)A$, where $A$
is a constant (matrix). Therefore we have $$L_{-d}L_d X_\nu=(\ol\nu+bd)(
\ol\nu+d-bd)AX_\nu.\eqno(3.40)$$ Since $x^{(1)}_0$ is an eigenvector of
$
L_{-d}L_d$, the first row of $A$ can only have a non-zero element in the
first position. Thus if we now choose $V^{(k,1)}={\sl Span}_\sC\langle
x^
{(1)}_\nu\ |\ \nu\in M_1\rangle$, then (3.40) tells us that $L_{-d}L_d
V^
{(k,1)}\subset V^{(k,1)}$. Now let's start from $k=k_0$ and choose
$V^{(k_
0,1)}$ as mentioned. By induction on $k>k_0$, we can choose $V^{(k,1)}$
to be the simple \HVira-submodule generated by $L_d V^{(k-1,1)}$. Using
$
L_{-d}L_d V^{(k_0,1)}\subset V^{(k_0,1)}$, by induction on $m$, we have
$
L_{md}V^{(k,1)}\subset V^{(k+m,1)}$ for $m,k\in\Z$ such that $k\ge
k_0,k+
m\ge k_0$. Since $V$ is a simple module, it must be generated by
$V^{(k_0,
1)}$. However, we have, by the arguments similar to the proof of Lemma
3.1, ${\sl dim\,}V_{a+\nu+kd}=1$ for $\nu\in M_1,k\ge k_0$. Thus by
Lemma
2.4, $V$ must be a module of the intermediate series. This completes the
proof of Theorem 2.2. \par\ \par\ \par \cl{\bf Reference} \ni\hi4ex\ha1
[1] V. Chari and A. Pressley, ``Unitary representations of the Virasoro
algebra and a conjecture of Kac,'' {\sl Compositio Math.}, {\bf
67}(1988), 315-342. \par\ni\hi4ex\ha1 [2] B. L. Feigin and D. B. Fuchs,
``Verma modules over the Virasoro algebra,'' in {\sl Lecture Notes in
Math.}, {\bf 1060}, pp.230-245, Springer, New York-Berlin, 1984. \par\ni
\hi4ex\ha1 [3] D. Friedan, Z. Qiu, and S. Shenker, ``Conformal
invariance, unitarity and two-dimensional critical exponents,'' in {\sl
Vertex Operators in Mathematics and Physics} (J. Lepowsky, S.
Mandelstam,
and I. M. Singer, Eds.), Springer, New York-Berlin, 1985.
\par\ni\hi4ex\ha
1 [4] ---, ``Details of the non-unitarity proof for highest weight
representations of the Virasoro algebra,'' {\sl Comm. Math. Phys.} {\bf
107}(1986), 535-542. \par\ni\hi4ex\ha1 [5] P. Goddard and D. Olive,
``Kac-Moody and Virasoro algebras in relation to Quantum physics,'' {\sl
Internat. J. Mod. Phys.}, {\bf A1}(1986), 303-414. \par\ni\hi4ex\ha1 [6]
V.G. Kac, ``Some problems on infinite-dimensional Lie algebras and their
representations,'' in {\sl Lie algebras and related topics, Lecture
Notes
in Math.} {\bf 933}, 117-126, Springer, New York-Berlin, 1982.
\par\ni\hi
4ex\ha1 [7] I. Kaplansky and L. J. Santharoubane, ``Harish-Chandra
modules over the Virasoro algebra,'' in {\sl Infinite-Dimensional Groups
with Application, Math. Sci. Res. Inst. Publ.} {\bf 4}(1985), 217-231.
\par\ni\hi4ex\ha1 [8] R. Langlands, ``On unitary representations of the
Virasoro algebra,'' in {\sl Infinite-Dimensional Lie algebras and Their
Applications} (S. N. Kass, Ed.), World Scientific, Singapore, 1986,
141-159. \par\ni\hi4ex\ha1 [9] C. Martin and A. Piard, ``Indecomposable
modules over the Virasoro Lie algebra and a conjecture of V. Kac,'' {\sl
Comm. Math. Phys.} {\bf 137}(1991), 109-132. \par\ni\hi4ex\ha1 [10] O.
Mathieu, ``Classification of Harish-Chandra modules over the Virasoro
Lie
algebra,'' {\sl Invent. Math.} {\bf 107}(1992), 225-234.
\par\ni\hi4ex\ha
1 [11] J. Patera \& H. Zassenhaus, ``The high rank Virasoro algebras,''
{\sl Comm. Math. Phys.} {\bf 136}(1991), 1-14. \par\ni\hi4ex\ha1 [12] C.
Rebbi, ``Dual models and relativistic quantum strings,'' {\sl Phys.
Rep.}
{\bf 12}(1974), 1-73. \par\ni\hi4ex\ha1 [13] Y. -C. Su, ``A
classification of indecomposable $sl_2(\C)$-modules and a conjecture of
Kac on Irreducible modules over the Virasoro algebra,'' {\sl J. Alg.},
{\bf 161}(1993), 33-46. \par\ni\hi4ex\ha1 [14] Y. -C. Su,
``Classification of Harish-Chandra modules over the super-Virasoro
algebras,'' {\sl Comm. Alg.} {\bf 23}(1995), 3653-3675.
\par\ni\hi4ex\ha1
[15] Y. -C. Su, ``Harish-Chandra modules of the intermediate series over
the high rank Virasoro algebras and high rank super-Virasoro algebras,''
{\sl J. Math. Phys.} {\bf 35}(1994), 2013-2023. \par\ni\hi4ex\ha1 [16]
Y.
-C. Su, ``Classification of indecomposable Harish-Chandra modules over
the Virasoro algebra and high rank Virasoro algebras,'' (to be
published)
\end{document}